\documentclass[12pt]{amsart}

\usepackage{amsmath,amssymb,esint,amsthm,bbm}
 \usepackage{xcolor}

\usepackage[colorlinks=true, pdfstartview=FitV, linkcolor=blue,
citecolor=blue, urlcolor=blue]{hyperref}

\calclayout

\newcommand{\M}{\mathcal{M}}

\newtheorem{thm}{\textbf{Theorem}}

\theoremstyle{remark}

\theoremstyle{definition}
\newtheoremstyle{Claim}{}{}{\itshape}{}{\itshape\bfseries}{:}{ }{#1}
\theoremstyle{Claim}

\title{A note on the limit of Orlicz norms} 

\author{David Cruz-Uribe, OFS, and Scott Rodney}

\address{David Cruz-Uribe, OFS \\
Dept. of Mathematics \\
University of Alabama \\
 Tuscaloosa, AL 35487, USA}
\email{dcruzuribe@ua.edu}

\address{Scott Rodney\\
Dept. of Mathematics, Physics and Geology \\ 
Cape Breton University \\
Sydney, NS B1Y3V3, CA} 
\email{scott\_rodney@cbu.ca}

\thanks{D.~Cruz-Uribe is supported by \
  research funds from the Dean of the College of Arts \& Sciences, the
  University of Alabama. S.~Rodney is supported by the NSERC Discovery
  Grant program. }

\date{November 30, 2020}



\begin{document}

\maketitle


A well-known result in measure theory is that if $(X,\M,\mu)$ is
a measure space, and if for some $r<\infty $,
\[ \|f\|_r  = \bigg(\int_X |f(x)|\,d\mu\bigg)^{\frac{1}{r}} < \infty, \]
then
\[ \lim_{p\rightarrow \infty} \|f\|_p = \|f\|_\infty. \]
(See Rudin~\cite[p.~73]{MR924157}.) The purpose of this
note is to generalize this result to the scale of Orlicz spaces.

To state our main result, we begin with some definitions and basic
facts about
Young functions and Orlicz spaces.  For complete information,
see~\cite{krasnoselskii-rutickii61, rao-ren91}; for a briefer summary,
see~\cite[Section~5.1]{MR2797562}.   Let $A :
[0,\infty) \rightarrow [0,\infty)$ be a Young function:  that is, $A$
is continuous, convex and increasing, $A(0)=0$, and $A(t)/t\rightarrow
\infty$ as $t\rightarrow \infty$.
Given a measurable, real-valued function on $f$, define the Orlicz
norm (more properly, the Luxemburg norm) by 
\[ \|f\|_A = \inf\bigg\{ \lambda > 0 :
  \int_X A\bigg(\frac{|f(x)|}{\lambda}\bigg)\,d\mu \leq 1 \bigg\}.  \]
When $A(t)=t^p$, $1\leq p<\infty$, then $\|f\|_A=\|f\|_p$.
Given two Young functions $A$ and $B$, we say $A\lesssim B$ if
there exists a constant $c$ such that for all $t>0$, $A(t)\leq
B(ct)$.  If $A\lesssim B$ and $B\lesssim A$, we write
$A\approx B$.    If $A\lesssim B$, then there exists a
constant $C>0$  depending only on $A$ and $B$, such that
$ \|f\|_A \leq C\|f\|_B$.

We are interested in the  Orlicz norms defined
with respect to the Young functions 
\[ B_{pq}(t) = t^p \log(e_0+t)^q, \]
where  $1\leq p < \infty$, $q>0$, and 
$e_0=e-1$.   Our result is the following.  

\begin{thm} \label{thm:main}
Given a measure space $(X,\M,\mu)$ and measurable function  $f$, if for some
$p\geq 1$ and $q_0>0$,  $\|f\|_{B_{pq_0}}<\infty$,  then
\begin{equation}\label{eqn:main1}
 \lim_{q\rightarrow \infty} \|f\|_{B_{pq}} = \|f\|_\infty.  
\end{equation}
\end{thm}

Before giving the proof, we make three remarks.
First, we emphasize that the significant feature of
Theorem~\ref{thm:main} is that the power $p$ remains
fixed, and it is only the logarithm term that grows.

Second, we first considered this result during an (unsuccessful)
attempt to generalize Moser iteration to the scale of Orlicz
spaces (see~\cite{DCU-SR20}, where we instead implemented De Giorgi
iteration in this setting).   For Moser iteration we needed to control the norm as $q\rightarrow \infty$.  We were
surprised that this result was not known, but a search of the literature did
not find it.

Third,
 in applications the Young functions $B_{pq}$ are often
defined with $e$ in place of $e_0$ (see, for
instance,~\cite{MR2797562}), and if we define
\[ \bar{B}_{pq}(t) = t^p \log(e+t)^q, \]
then $B\approx \bar{B}$.  However, in order to get an equality in the
conclusion of Theorem~\ref{thm:main} we are required to assume that
$B_{pq}(1)=1$.

\medskip

\begin{proof}
  We prove \eqref{eqn:main1} in two steps.  We will first show that
  \begin{equation} \label{eqn:liminf}
\liminf_{q\rightarrow \infty} \|f\|_{B_{pq}} \geq
    \|f\|_\infty.  
  \end{equation}
  If $\|f\|_\infty=0$, there is nothing to prove.  Now suppose that
  $0<\|f\|_\infty<\infty$.  Since for some $q_0>0$, $\|f\|_{B_{pq_0}}<\infty$, the
  level sets of $f$, $S(\lambda) = \{ x\in X : |f(x)|>\lambda\}$,
  $\lambda>0$, have finite measure: this follows from the
  corresponding result for $L^p$ spaces, since
  $\|f\|_p \lesssim \|f\|_{B_{pq_0}}$.  Fix $\epsilon>0$ and let
  $M= S(\|f\|_\infty-\epsilon)$.  Then $0<\mu(M)<\infty$, and so
  \[ \|f\|_{B_{pq}} \geq \|f\chi_M\|_{B_{pq}} \geq
    (\|f\|_\infty-\epsilon) \|\chi_M\|_{B_{pq}}. \]
 Since $\epsilon>0$ is arbitrary,  to prove
 inequality~\eqref{eqn:liminf} it will suffice to show 
  \begin{equation} \label{eqn:liminf2}
    \liminf_{q\rightarrow \infty} \|\chi_M\|_{B_{pq}} \geq 1. 
  \end{equation}

 To estimate $\|\chi_M\|_{B_{pq}}$, fix $q>q_0$.   By the
  definition of the Orlicz norm,
  \[ \lambda = \|\chi_M\|_{B_{pq}} = B_{pq}^{-1}(\mu(M)^{-1})^{-1}.  \]
  Hence,
  \[  \mu(M)^{-1} = B_{pq}(\lambda^{-1}) =
    \lambda^{-p}\log(e_0+\lambda^{-1})^q. \]
  If  we always have that $\lambda \geq 1$, then \eqref{eqn:liminf2}
  clearly holds.  If for some $q$, $\lambda<1$,  fix $\delta>0$ such
  that  $\log(e_0+\lambda^{-1})=1+\delta$;  by the binomial theorem,
  \[ \lambda^{-p}\log(e_0+\lambda^{-1})^q  = (e^{1+\delta}-e_0)^p (1+\delta)^q \geq  1+q\delta > q\delta.  \]
  Thus, $\delta < \frac{\mu(M)^{-1}}{q}$.
  If we solve for  $\lambda$ in the definition of $\delta$, we get
  \[  \|\chi_M\|_{B_{pq}}  = \lambda > \big[
    \exp\big(1+\tfrac{\mu(M)^{-1}}{q}\big)-e_0 \big]^{-1}.  \]
Therefore,   if we take the limit infimum as $q\rightarrow\infty$, we get \eqref{eqn:liminf2}.

Finally, suppose $f$ is unbounded.  For   any $N>1$ define
  $f_N=\min(|f|,N)$. Then the above argument shows that
  \[ \liminf_{q\rightarrow \infty} \|f\|_{B_{pq}}
    \geq \liminf_{q\rightarrow \infty} \|f_N\|_{B_{pq}} \geq
    \|f_N\|_\infty = N,  \]
and  if we let $N\rightarrow \infty$, \eqref{eqn:liminf} follows. 

  \medskip

 We will now show that
   \begin{equation} \label{eqn:limsup}
\limsup_{q\rightarrow \infty} \|f\|_{B_{pq}} \leq
    \|f\|_\infty.  
  \end{equation}
  We may assume that $0<\|f\|_\infty<\infty$; otherwise
  \eqref{eqn:limsup} is immediate.  Since $\|f\|_{B_{pq_0}}<\infty$, by the
  definition of the norm we have that there exists $\lambda>0$ such that
  \[ \big(\tfrac{|f|}{\lambda}\big)^p\log\big(e_0+
    \tfrac{|f|}{\lambda}\big)^{q_0} \in L^1(X).  \]
Moreover,  this is true for any $\lambda$,
  $0<\lambda<\infty$; this follows at once from the fact that for any
  $c>0$, the function
  \[ \frac{\log(e_0+t)}{\log(e_0+ct)} \]
  is bounded and bounded away from $0$ for all $t>0$.

  In particular, if we fix $\epsilon>0$ and let
  $\lambda=(1+\epsilon)\|f\|_\infty$, then for almost every $x\in X$,
  and for all $q>q_0$,  since $e_0+(1+\epsilon)^{-1}<e$,
  \[ B_{pq}\big(\tfrac{|f(x)|}{\lambda}\big) \leq
    B_{pq_0}\big(\tfrac{|f(x)|}{\lambda}\big).  \]
  Furthermore,
  \[ B_{pq}\big(\tfrac{|f(x)|}{\lambda}\big)
    \leq (1+\epsilon)^{-p} \log(e_0+(1+\epsilon)^{-1})^q.  \]
  Again since $e_0+(1+\epsilon)^{-1}<e$, the right-hand side tends to $0$ as
  $q\rightarrow \infty$.  Hence, by the dominated convergence
  theorem, for all $q>q_0$ sufficiently large,
  \[ \int_X B_{pq}\bigg(\frac{|f(x)|}{\lambda}\bigg)\,d\mu \leq 1. \]
Therefore, by the definition of the Orlicz norm, for all such $q$,
  \[ \|f\|_{B_{pq}} \leq (1+\epsilon)\|f\|_\infty, \]
and   inequality~\eqref{eqn:limsup} follows at once.   This completes
the proof.

\end{proof}

\bibliographystyle{plain}
\bibliography{Limit-Orlicz}

\end{document}